\documentclass[10pt]{amsart}
\usepackage{amsmath,amssymb}

\newcommand\inv{^{-1}}

\newcommand\fh{\mathfrak h}\newcommand\fj{\mathfrak j}

\newcommand\fz{\mathfrak z}

\newcommand\fm{\mathfrak m}
\newcommand\ff{\mathfrak F}

\newcommand\G{\Gamma}
\newcommand\g{\gamma}

\newcommand\Ra{\mathbb R}

\newcommand\Za{\mathbb Z}

\newcommand\fd{\mathfrak d}

 \DeclareMathOperator{\Id}{Id}
 \DeclareMathOperator{\diag}{diag}
\DeclareMathOperator{\Diff}{Diff} 
\DeclareMathOperator{\Aff}{Aff}

\DeclareMathOperator{\Hom}{Hom} \DeclareMathOperator{\Aut}{Aut}
\DeclareMathOperator{\Ad}{Ad}

\newtheorem{theorem}{Theorem}[section]
\newtheorem{proposition}[theorem]{Proposition}
\newtheorem{lemma}[theorem]{Lemma}
\newtheorem{corollary}[theorem]{Corollary}

\begin{document}

\title{Deformations of group actions}
\author{David Fisher}
\thanks{Author partially supported by NSF grant DMS-0226121 and a PSC-CUNY grant.}

\begin{abstract}
Let $G$ be a noncompact real algebraic group and $\G<G$ a lattice.
One purpose of this paper is to show that there is a smooth, volume
preserving, mixing action of $G$ or $\G$ on a compact manifold which
admits a smooth deformation.  In fact, we prove a stronger statement
by exhibiting large finite dimensional spaces of deformations. We
also describe some other, rather special, deformations when
$G=SO(1,n)$.
\end{abstract}

\maketitle

\section{Introduction}

In recent years, many results have been proven concerning local
rigidity of group actions. Most known results concern actions of
semi-simple Lie groups  with all simple factors of real rank at
least two and lattices in such groups. For the best known results
in that category and some historical references, see \cite{FM1}.
More recently some partial results have been proven for more
general groups, including lattices in $Sp(1,n)$ and $F_4^{-20}$,
see \cite{FM2} and \cite{H}.  In this context it is interesting to
ask when Lie groups and their lattices have actions by
diffeomorphisms which admit non-trivial perturbations or
deformations.  Actions of $SL(n,\Za)$  constructed by Katok and
Lewis admit non-trivial deformations and their construction was
modified by Benveniste to produce actions of some higher rank
simple Lie groups with non-trivial volume preserving deformations
\cite{B,KL}.  Benveniste's deformations are also non-trivial when
all actions in question are restricted to any lattice in the
acting group.  Throughout this paper, by a connected real
algebraic group, I mean the connected component of the real points
of an algebraic group defined over $\Ra$.  One of the main results
of this paper is the following generalization of Benveniste's
theorem:

\begin{theorem}
\label{theorem:simple} Let $G$ be a non-compact, simple, connected
real algebraic group and $\G<G$ a lattice. Then there is an a
smooth, volume preserving, mixing, action of $G$ which admits
non-trivial, volume preserving deformations. Furthermore, the
deformations remain non-trivial when restricted to $\G$.
\end{theorem}

{\bf Remarks:}
\begin{enumerate}
\item It is immediate from the construction that we can construct
actions where the space of deformations has arbitrarily large,
finite, dimension. It would be interesting to know when the
deformation space is infinite dimensional. For any group $G$ which
contains a lattice with a homomorphism onto $\Za$, it is easy to
build examples with infinite dimensional families of perturbation
using simple induction constructions and Lemma \ref{lemma:induction}
below.

\item The construction of the actions and deformations is a
refinement of Benveniste's which is based on an earlier
construction of Katok and Lewis \cite{B,KL}. The proof that the
deformations are non-trivial is new and uses a consequence of
Ratner's theorem due to N.Shah. Benveniste's proof used
commutativity coming from a higher rank assumption on $G$.

\item Benveniste's explicit deformations are non-trivial in a
sense that is weaker than what is meant by Theorem
\ref{theorem:simple}, see below for discussion.
\end{enumerate}

We now introduce some basic notions in order to describe in detail
the meaning of the statement of Theorem \ref{theorem:simple} and
also to clarify Remark $(1)$ and show that we produce ``many"
deformations.  Let $\rho$ be an action of a group $D$ on a manifold
$M$. Then the space of actions of $D$ on $M$ is naturally the space
$\Hom(D,\Diff(M))$. The group $\Diff(M)$ acts on this space by
conjugation, and two group actions are conjugate if and only if they
are in the same $\Diff(M)$ orbit.  Note that if $D$ is finitely
generated, $\Hom(D,\Diff(M))$ is naturally a closed subset of the
product of a finite number of copies of $\Diff(M)$. In the finite
dimensional setting, i.e. where we replace $\Diff(M)$ by an
algebraic group $H$, the space $\Hom(D,H)$ is an algebraic variety,
and there is a construction of a quotient scheme $\Hom(D,H)//H$
which allows one to study the space of non-trivial deformations.  In
the infinite dimensional setting it is unclear that there is any
meaningful algebraic geometry, or that $\Hom(D,\Diff(M))$ is a
manifold at ``most" (or even any) points. We now describe in more
detail the space of deformations produced in the proof of Theorem
\ref{theorem:simple}.  We do this for $\G$ actions, where the space
$\Hom(\G,\Diff(M))$ is at least a closed subset of a comprehensible
topological space, the $k$-fold product
$\Diff(M){\times}{\cdots}{\times}{\Diff(M)}$ where $k$ is the number
of generators of $\G$. (There are various ways in which one can
topologize $\Hom(G,\Diff(M))$ but this is more complicated, and we
will not discuss it here.) Our construction produces an action of
$\G$ on a manifold $M$ and a collection of deformations which
provide an embedding of a small, finite dimensional ball in
$\Hom(\G,\Diff(M))$ which intersects the $\Diff(M)$ orbit of $\rho$
only at $\rho$. Benveniste's construction, even using our proof in
place of his, only produces an action $\rho$ with deformations which
define a map of a ball into $\Hom(\G,\Diff(M))$ which intersects the
$\Diff(M)$ orbit of $\rho$ at a countable set which could, a priori,
accumulate on $\rho$. Since we do not know, in this context, that
the representation variety is even locally a manifold, it is unclear
that the two statements are equivalent. As will be clear below,
given $\G$ and a positive integer $n$ we can construct a manifold
$M$ and an action $\rho$ such that the ball discussed above has
dimension $n$.  By ``blowing up and gluing" along many distinct
pairs of closed orbits, one can even arrange that $n$ is much larger
than $\dim(M)$.

In addition to removing the rank restriction in Benveniste's
theorem, our use of Ratner theory is sufficiently robust to allow
us to use slight modifications of the examples for Theorem
\ref{theorem:simple} to prove the following stronger theorem:

\begin{theorem}
\label{theorem:general} Let $G$ be a connected, non-compact real
algebraic group and $\G<G$ a lattice. Then there is a faithful,
smooth, volume preserving, ergodic, action of $G$ which admits
non-trivial, volume preserving deformations. Furthermore, the
deformations remain non-trivial when restricted to $\G$.
\end{theorem}

\noindent {\bf Remarks:} \begin{enumerate} \item As will be clear
in the proof, Theorem \ref{theorem:general} is only really
difficult if the abelianization of $G$ is compact. See Proposition
\ref{proposition:suffices}. \item Theorem \ref{theorem:general} is
in a sense sharp, since any smooth action of a compact group is
locally rigid by a theorem of Palais \cite{P}. \item Though the
author knows of no claims in the literature, it is quite likely
that one can prove local rigidity theorems for groups which are
not semisimple. In particular, the results of \cite{FM1} seem
likely to hold for $G$ any algebraic group all of whose
non-trivial factors are semisimple groups with all simple factors
of real rank at least two and for lattices $\G$ in such $G$. Any
such $G$ or $\G$ can be shown to have property $(T)$, so the
arguments in \cite{FM2} go through verbatim.  The main difficulty
in proving local rigidity for quasi-affine actions appears to be
in generalizing some of the results of \cite{FM0}. The results in
\cite{W2} should be relevant to this generalization.
\end{enumerate}

With the exception of some of the results in \cite{Be2,FM1,FM2},
all local rigidity theorems known to the author are for actions on
manifolds which are homogeneous.  The deformations constructed to
prove Theorems \ref{theorem:simple} and \ref{theorem:general}
above are on manifolds which are not homogeneous.  In section
\ref{section:extradeformations} we show that some affine actions
of lattices in $SO(1,n)$ on homogeneous manifolds have non-trivial
deformations.  This then implies that the same is true for some
affine actions of $SO(1,n)$ as well.

One motivation for this paper is the so-called Zimmer program to
classify actions of higher rank lattices, see
\cite{Zim_ICM,Zim_Most} for more details.  In the context of both
this program and other older rigidity results, it is not
surprising that actions of rank $1$ lattices exhibit greater
flexibility than higher rank lattices.

\noindent
 {\it Acknowledgements:} I would like to thank
E.J.Benveniste and K.Whyte for many useful conversations
concerning the construction and properties of ``exotic actions"
discussed in this paper.  I learned a great deal about these
examples during the work on \cite{BF} and \cite{FW}. I would also
like to thank Dave Witte Morris for some useful emails concerning
algebraic groups and lattices. A proof that ``bending
deformations" gave non-trivial deformations of affine actions of
lattices and Lie groups was discovered jointly with R.Spatzier
during an extremely pleasant visit to the University of Michigan
in Spring of 2002. This proof uses the normal form theory of
Hirsch-Pugh-Shub for partially hyperbolic diffeomorphisms and is
sketched in section 5. That proof led me to the proof using Ratner
theory given in section 5, which was the inspiration for the use
of Ratner theory in the proofs of Theorem \ref{theorem:simple} and
Theorem \ref{theorem:general}. Theorem \ref{theorem:simple} was
motivated by a question of Y.Shalom concerning the existence of
deformations of actions of lattices in $SU(1,n)$ and $SP(1,n)$.

%Finally, the proof in section \ref{section:lrforcompact} is more
%or less the proof of \cite[Theorem 1.1]{FM2} given in Section
%$5.1$ of that paper. I include a proof here merely because the
%argument simplifies dramatically for compact groups. I would like
%to take this opportunity to thank G.A.Margulis for the pleasure of
%the collaboration that produced the papers \cite{FM1,FM2,FM0}.

\section{Preliminaries}
\label{section:preliminaries}

\subsection{Invariant measures and centralizers of group actions}
\label{subsection:ratner}

In this section we recall a consequence of work of Ratner and Shah
on invariant measures which is instrumental in our proofs. For the
proofs of Theorems \ref{theorem:simple} and \ref{theorem:general},
our use of this result will be based on the fact that the
perturbations we construct are conjugate back to the original
action on a set of full measure.  Other applications in Section
\ref{section:extradeformations} will be to deformations that are
trivial by construction on a ``large" subgroup of the acting
group.  We let $\Aff(H/{\Lambda})$ denote the affine group of
$H/{\Lambda}$.   This group consists of diffeomorphisms of
$H/{\Lambda}$ which are covered by maps of the form $A{\circ}L_h$
where $A$ is an automorphism of $H$ such that $A(\Lambda)=\Lambda$
and $L_h$ is left translation by $H$.  The group
$\Aff(H/{\Lambda})$ is a Lie group and is, in fact a quotient of a
a subgroup of $\Aut(H){\ltimes}H$, see \cite{FM0} for more
discussion.

\begin{corollary}
\label{corollary:rsw} Let $H$ be a connected real algebraic group,
$\Lambda<H$ a lattice, $\mu_{H}$ the measure on $H/{\Lambda}$
induced by a fixed Haar measure on $H$ and $F<H$ a group that
either
\begin{enumerate}
\item contains a subgroup $F'<F$ generated by unipotent elements,
\item contains a subgroup $F'<F$ such that $F'$ is a lattice in a
subgroup $F''<H$ generated by unipotent elements.
\end{enumerate}
Further assume $F'$ acts ergodically on $H/{\Lambda}$. If $\phi$
is an essentially surjective, essentially injective measurable map
from $H/{\Lambda}$ to $H/{\Lambda}$ which commutes with the $F$
action, then $\phi$ is translation by an element of the
centralizer of $F'$ in $\Aff(H/{\Lambda})$.
\end{corollary}

\noindent Though this corollary is essentially noted in Witte's
paper \cite{W} where a more general result is proven.  For
completeness, we briefly recall the proof.

\begin{proof}
If $\phi$ commutes with $F'$, then the graph $N$ of $\phi$ is
preserved by $F'$. Letting $\tilde \phi:H/\Lambda{\rightarrow}
H/{\Lambda}{\times}H/{\Lambda}$ be given by $\tilde
\phi=(\Id,\phi)$, it is straightforward to check that $\tilde
\phi_*\mu_H$ is an $F'$ invariant, ergodic measure on $N$ that
projects to $\mu_H$ on each $H/{\Lambda}$.  Given a subgroup $L<H$,
we denote by $\Delta(L)$ the diagonal embedding of $L$ in
$H{\times}H$.
%Furthermore, since $F'$ contains elements of $H$
%which generate unbounded subgroups, the $F'$ action on
%$H/{\Lambda}$ is ergodic and so is the $F'$ action on $(N,\tilde
%z_*\mu_H)$.
By Shah's extensions of Ratner's theorem to groups generated by
unipotents and their lattices \cite{R,S},  this implies that that
$N$ is a closed $L$ orbit in $H/{\Lambda}{\times}H/{\Lambda}$ where
$L$ is a closed subgroup of $H{\times}H$ containing $\Delta(F')$.
Since $\tilde \phi$ composed with projection on the second factor is
essentially injective and essentially surjective, we have that $L$
is exactly $\Delta(H)$. This implies that $N$ is a closed
$\Delta(H)$ orbit in $H/{\Lambda}{\times}H/{\Lambda}$, which implies
that $\phi$ is translation by an element of $H$ composed with an
automorphism of $H$. This implies that $\phi$ is in fact an element
of $Z_{\Aff(H/{\Lambda})}(F')$.
\end{proof}

\subsection{Reductions for the proof of Theorem
\ref{theorem:general}} \label{subsection:reductions}

We recall briefly some facts concerning (real) algebraic groups.
Given a connected real algebraic group $G$, we can write $G$ as a
semi-direct product $L{\ltimes}U$ where $U,L<G$ and $U$ is the
unipotent radical of $G$ and $L$ is a reductive real algebraic
group which is called a Levi complement.

\begin{proposition}
\label{proposition:suffices}
To prove Theorem \ref{theorem:general}, it suffices to consider the case
when the Levi complement of $G$ has compact center.
\end{proposition}

\begin{proof}
We assume the center of the Levi complement is non-compact and
produce the family of perturbations required in Theorem
\ref{theorem:general}. We have a projection
$\pi_1:G{\rightarrow}L$, so any $L$ action defines a $G$ action.
We can write $L$ as $SZ(L)$ where $S$ is semisimple,
$S{\cap}Z(L)=F$ is finite and the product is almost direct.
Therefore $S$ is normal in $L$ and there is a projection
$\pi_2:L{\rightarrow}Z(L)/F$.  The group $Z(L)/F$ is an abelian
Lie group and non-compact whenever $Z(L)$ is non-compact, so we
can find a third projection $\pi_3:Z(L)/F{\rightarrow}\Ra$.  We
write $\pi=\pi_3{\circ}\pi_2{\circ}\pi_1$ and note that any $\Ra$
action defines a $G$ action via composition with $\pi$.

There are many $\Ra$ actions which are smooth and mixing and which have non-trivial perturbations,
we give an example for completeness. Let $\Lambda<SL(2,\Ra)$ be a cocompact lattice and let $\rho_{\Ra}$
 be the action of $\Ra$ on
$SL(2,\Ra)/\Lambda$ defined by the horocycle flow.  (I.e. defined
by identifying $\Ra$ with strictly upper triangular matrices in
$SL(2,\Ra)$ and acting by left translation.)  It is easy to
perturb this action, even as an action by left translations, since
this action has zero entropy and many nearby actions defined by
``nearby" subgroups of $SL(2,\Ra)$ have positive entropy. It is
straightforward to construct families on which the entropy is a
strictly increasing as one moves away from $\rho_{\Ra}$.  Since
Lyapunov exponents are also conjugacy invariants of
diffeomorphisms, by choosing an $\Ra$ action on a larger manifold,
one can construct large dimensional families of deformations which
are easily seen to be non-trivial and to give rise to non-trivial
deformations of the faithful actions described in the next
paragraph.

To construct a faithful action of $G$ with pertubations, we embed
$G$ in $SL(N+1,\Ra)$, choose a cocompact lattice
$\Delta<SL(N+1,\Ra)$ and let $\rho_F$ be the left translation
action of $G$ on $SL(N+1,\Ra)/{\Delta}$.  We then form
$M=SL(2,\Ra)/\Lambda{\times}SL(N+1,\Ra)/{\Delta}$ and let $\rho$
be the diagonal $G$ action defined by $\pi{\circ}{\rho_{\Ra}}$ and
$\rho_F$. It is straightforward to check that this satisfies the
hypotheses of Theorem \ref{theorem:general}. The fact that the
action is mixing follows from the Howe-Moore theorem, see e.g.
\cite{ZB}. It is also straightforward to check that $\rho|_{\G}$
provides the necessary action of $\G$ when $\G<G$ is a lattice.
\end{proof}

We now turn to the case where the entire Levi complement is compact.
In order to prove the existence of deformations in this case, we
will need to use induced actions.  Let $F$ be a locally compact
topological group and $D<F$ a closed subgroup.  Assume $D$ acts
continuously on a compact space $X$.  We can then form the space
$(F{\times}X)/D$ where $D$ acts by $(f,x)d=(fd{\inv},dx)$. This
space is clearly equipped with a left $F$ action which we refer to
as the {\em induced $F$ action}.  It is easy to check that if $F$ is
a Lie group, $X$ is a smooth manifold  and $D$ acts smoothly on $X$,
then the induced $F$ action is also smooth.

\begin{lemma}
\label{lemma:induction} If $H$ is a real Lie group and $L<H$ is a
closed cocompact subgroup, then if $L$ has an action on a compact
manifold $M$ which admits non-trivial deformations (resp.
perturbations), the induced $H$ action on $(H{\times}M)/{L}$ also
admits non-trivial deformations (resp. perturbations).
\end{lemma}

\begin{proof}
Let $\rho_0$ be an $L$ action on $M$ and $\rho_t$ a non-trivial one
parameter deformation.  We can form the induced actions $\tilde
\rho_t$ of $H$ on $(H{\times}M)/{L}$ which are again smooth actions
and it is clear that $\tilde \rho_t$ is a one parameter deformation
of $\tilde \rho_0$.  For a contradiction, assume that deformations
$\tilde \rho_t$ are trivial.  Then there exist diffeomorphisms
$\phi_t$ of $(H{\times}M)/{L}$ such that
$\phi_t{\circ}\rho_0=\rho_t{\circ}\phi_t$.  Let
$\pi:(H{\times}M)/{L}{\rightarrow}H/{L}$ be the natural projection.
Then $\phi_t{\circ}\pi$ is an $H$ equivariant map from
$((H{\times}M)/{L},\rho_t)$ to $(H/{L}, \rho_H)$ which is close to
$\pi$ as a smooth map.  The sets $(\phi_t{\circ}\pi){\inv}([h])$ are
just the sets $\phi_t\inv((\pi{\inv}[h]))$ and so are an $H$
invariant foliation of $(H{\times}M)/{L}$.  In particular
$\phi_t\inv((\pi{\inv}[e]))$ is a $L$ invariant copy of $M$, and it
is easy to check that the $L$ action on $\phi_t{\inv}(\pi{\inv}[e])$
is exactly $\rho_t$. Therefore $\phi_t|_{\pi{\inv}([e])}$ is a
conjugacy between $\rho_0$ and $\rho_t$.  It is straightforward to
check that $\phi_t|_{\pi{\inv}([e])}$  is indeed a small, smooth
path in $\Diff(M)$.  The same proof applies to perturbations in
place of deformations.
\end{proof}

\begin{proposition}
\label{proposition:noncompact}
It suffices to prove Theorem \ref{theorem:general} when the Levi complement of $G$ is
not compact.
\end{proposition}

\begin{proof} We prove that if the Levi complement is compact, then
there are deformations. By this assumption, $G=C{\ltimes}U$ where
$C$ is compact and $U$ is unipotent.   There exist non-trivial
homomorphisms $\sigma:U{\rightarrow}\Ra$.  As in the proof of
Proposition \ref{proposition:suffices}, we can use an $\Ra$ action
$\rho_{\Ra}$ with non-trivial perturbations to construct a $U$
action with non-trivial perturbations.  This yields a $G$ action
with non-trivial perturbations by Lemma \ref{lemma:induction}. The
$G$ action is volume preserving if $\rho_{\Ra}$ is, since there is
a $G$ invariant volume on $G/U$.  It is easy to check that this
action is mixing when $\rho_{\Ra}$ is mixing.  To obtain a
faithful, mixing action, one proceeds as in Proposition
\ref{proposition:suffices}.  To verify that the deformations are
non-trivial on the faithful action requires some care, and we
sketch a simple argument using entropy or Lyapunov exponents. Let
$B<U$ be a one parameter subgroup mapped onto $\Ra$ by $\sigma$.
Then the restriction of the induced action to $B$ is the product
of the trivial $B$ action on $G/U$ with the action
$\rho_{\Ra}{\circ}{\sigma}|_B$.  This shows that the Lyapunov
exponents or entropy of $\rho_{\Ra}$ are an invariant of the
induced action, and therefore of the diagonal action used to
define a faithful $G$ action.

Given $\G<G$, a theorem of Auslander shows that $\G{\cap}U$ is a
lattice in $U$ and therefore that $\sigma(\G{\cap}U)$ is infinite
and unbounded \cite{A}.  From this it is easy to check that the
restriction of the above action to $\G$ proves Theorem
\ref{theorem:general} for $\G$ actions.
\end{proof}

\noindent{\bf Remarks:}
\begin{enumerate}
\item The techniques of this subsection are really only necessary
for the case of a noncompact real algebraic group all of whose
simple quotients are compact.  As long as there is a non-compact
simple quotient, one can proceed by the methods of section
\ref{section:deformations}.

\item Here we only use rough dynamical invariants such as entropy
and Lyapunov exponents to see that perturbations are non-trivial.
Using finer dynamical invariants, one can easily see that the
spaces of deformations yielded by Proposition
\ref{proposition:suffices} and Proposition
\ref{proposition:noncompact} are infinite dimensional. We do not
deal with this here, since it would lead to a long digression on
classical dynamics.
\end{enumerate}

\section{Constructing ``exotic actions"}
 \label{section:jerome'sconstruction}

We work with a variant of a construction developed by Benveniste in
\cite{B}, which is inspired by the construction in \cite{KL}. Let
$H=SL(n,\Ra)$ and $J=SL(m,\Ra)$ with $m<n-1$. View $J<H$ as a
subgroup via the standard embedding as block $m$ by $m$ matrices in
the upper left corner. For applications we will choose $m$ large
enough so $G<SL(m,\Ra)$=J. We let $Z=SL(n-m,\Ra)<Z_H(J)$ be the
block $n-m$ by $n-m$ matrices in the lower right hand corner. In
fact $Z_H(J)=Z(H)DZ$ where $Z(H)$ is the center of $H$ and $D$ is
the group of diagonal matrices in $H$ commuting with $J{\times}Z$,
i.e. the group of diagonal matrices of the form
$\diag\{d_1,\ldots,d_m,d_{m+1},ldots,d_{n}\}$ with
$d_1=d_2=\cdots=d_m$ and $d_{m+1}={\cdots}=d_{n}$.  Note that $D<H$
forces $d_1^{m}d_{n+1}^{n-m}=1$.  We recall the existence of certain
types of lattices in $H$.

\begin{lemma}
\label{lemma:latticeexists} If $n-m$ is a multiple of $m$, there
exists a cocompact lattice $\Lambda<H$ such that:

\begin{enumerate}
\item $\Lambda{\cap}J=\Lambda_J$ is a cocompact lattice in $J$
\item $\Lambda{\cap}Z=\Lambda_Z$ is a cocompact lattice in $Z$
\item $\Lambda{\cap}D$ is trivial
\item $\Lambda{\cap}DJZ=\Lambda_J {\times} \Lambda_Z$
\end{enumerate}

\end{lemma}

\begin{proof}
This follows easily from general methods to construct cocompact
lattices using unitary groups of Hermitian forms.  See for example
\cite[Proposition 7.54 and Figure 10.2]{Mo} for a construction that
easily yields lattices where the first two items are clear.  Since
$D$ is diagonalizable over $\mathbb Q$ (i.e. $\mathbb Q$ isotropic),
point $(3)$ follows from \cite[Theorem $4.11$]{PR}.  Point $4$
follows from \cite[Theorem 4.13]{PR} or rather from part of it's
proof and a use of change of base to realize $\Lambda$ as the
integer points for a $\mathbb Q$-structure on $H$ that restricts to
a $\mathbb Q$-structure on $DJZ$.
\end{proof}

Note that $JZ=SL(m,\Ra){\times}SL(n-m,\Ra)<SL(n,\Ra)$ is a closed
subgroup such that $JZ{\cap}{\Lambda}=\Lambda_{JZ}$ is a lattice in
$JZ$. To simplify some arguments, we pass to a subgroup of finite
index in $\Lambda$ to guarantee that $\Lambda$ does not intersect
the center of $J$,$H$ or $Z$.  Let $M = H/ \Lambda$. If $h \in H$,
we will write $[h]$ for the image of $h$ in $M$.  The following is
now easy from our choices of groups, subgroups and lattice:

\begin{corollary}
\label{corollary:orbits} The orbit $Z(H)JZD[e]$ is one-to-one
immersed.  The orbit $J[e]$ is an embedded submanifold $J$
equivariantly diffeomorphic to $J/{\Lambda_J}$ and the orbit
$JZ[e]$ is an embedded submanifold which is $JZ$ equivariantly
diffeomorphic to $J/{\Lambda_J}{\times}Z/{\Lambda_Z}$. In
particular, if $z{\in}Z_H(J)$ and $zJ[e] \subset JZ[e]$ then
$z{\in}Z$.
\end{corollary}

\begin{proof}
It is a standard fact that the orbit of a point $x$ in a
homogeneous space for a Lie group $H$ under a closed subgroup $L$
is an immersed submanifold which is equivariantly diffeomorphic to
$L/L_x$ where $L_x$ is the stabilizer of $x$ in $L$.  Furthermore,
the orbit is embedded if $L/L_x$ is compact. Combined with Lemma
\ref{lemma:latticeexists}, this yields the corollary immediately.
\end{proof}

\noindent{\bf Remark:} The construction in \cite{B} involves less
careful choices of subgroups and lattices and this care is taken
here exactly to allow us to show the existence of an embedded disk
in the space of deformations.

Note that for $z \in Z$ small enough, $zJ[e]$ does not intersect
the closed set $J[e]$. Let $C_0 \subset H/ \Lambda$ be the orbit
$J[e]$ and $C_1 \subset H/ \Lambda$ be the orbit $J[z]$, where the
action of $J$ is the obvious left action on $H / \Lambda$; the
sets $C_0$ and $C_1$ are closed submanifolds (by our assumption on
$\Lambda$) and are disjoint. We now form a manifold $Y_z$ by
blowing up $C_0$ and $C_1$ and gluing the resulting exceptional
divisors. More precisely, define $Y_z$ as follows.

We may identify $TM$ with $M \times \frak{h}$, and then the
derivative action of $H$ on $TM$ is given by $h(m, X) = (hm, Ad(h)
X)$. The tangent bundle to the $J$-orbits is then clearly $M
\times \frak{j}$; since the adjoint action of $J$ on $\frak{h}$ is
reductive, there is an $J$-invariant complement $V$ to $\frak{j}$,
and $M \times V$ is an $J$- invariant subbundle of $TM$ which
restricts to a normal bundle for $C_0$ or $C_1$.  Define maps

\[
\displaystyle \phi_i : J \times _{\Lambda_J} V \rightarrow M, i
=0,1
\]

by

\[
\displaystyle \phi_0 ([j], v) = j\exp(v)[e], \phi_1 (h,v) = j\exp(
v)[z];
\]
These maps are clearly well-defined and $J$-equivariant, and there
are open neighborhoods $U_i \subset J \times_{\Lambda_J} V$ of $J
\times _{\Lambda_J} 0$ such that $\phi_i$ restricted to $U_i$ is a
diffeomorphism onto a neighborhood of $C_i$, for $i=0,1$. Now let
$S^+ = (V -{0}) / {\Ra}^+$ and let $L ^+ \rightarrow S^+$ be the
tautological bundle of (closed) rays over $S^+$, i.e.,

\[
\displaystyle L^+ = \{(v, [w]) \in V \times S ^+ \vert v = rw, r \geq
0 \}.
\]
Denote the zero-section of $L^+$ by $S$. Also, let $L \rightarrow
S^+$ denote the tautological line bundle; observe that $L$ is two
copies of $L^+$ glued along $S$. The group $GL(V)$ acts on $L,
L^+$, $S^+$, and $S$ and therefore $J$ does as well. We can then
form the space $J \times _{\Lambda_J} L^+$, and there is an
$J$-equivariant embedding $V-0 \to L ^+$. Let $\hat {U}_i = (U_i -
(J \times _{\Lambda_J} 0)) \cup S$, where $U_i$ is regarded as a
subset of $J \times _{\Lambda_J} L^+$ through the above embedding
for $i=0,1$. Now let

\[
\displaystyle M_z = (M - C_0 \cup C_1) \sqcup \hat{U}_0 \sqcup
\hat{U}_1) / \mathcal{R}
\]
 where the equivalence relation ${\mathcal R}$ is generated by the
requiring
 $u \equiv \phi_i (u)$
for $ u \in U_i$ and $\phi_i(u) \in M - (C_0 \cup C_1)$.  $M_z$ is
then a manifold with boundary with two boundary components $B_0$
and $B_1$, each diffeomorphic to $J \times _{\Lambda_J} S^+$; it
is also equipped with an $J$-action by declaring the action to be
the restriction of the action on  $M$ to $M_z - (B_0 \cup B_1) = M
- (C_0 \cup C_1)$, and the obvious action on $B_i = J \times
_{\Lambda_J} S^+$ for $i= 0, 1$. This action is smooth by
construction.

Now let $Y_z$ be the manifold obtained by gluing $B_0$ to $B_1$ in
the following way: identifying each of these boundary components
to $J \times _{\Lambda_J} S^+$, glue the point $(j,s)$ on $B_0$ to
the point $(j, -s)$ on $B_1$; this is obviously well-defined. It
is obvious too that the action on $M_z$ descends to $Y_z$. Let
$B_z$ be  the image of the $B_i$s in $Y_z$, and let $U$ be $Y_z -
B = M - (C_0 \cup C_1)$; $B$ is a closed $J$-invariant set and $U$
is an open dense $J$-invariant set. Note that there is a $J$-map
$\phi : J \times _{\Lambda_J} L \rightarrow Y _z$
 and neighborhoods $W$ of $J \times _{\Lambda_J} S ^+$ in $J \times
_{\Lambda_J} L$ and $W ^{\prime}$ of $B$ in $Y _z$ such that
 $\phi \vert _W : W \rightarrow W ^{\prime}$
is a diffeomorphism.

The Haar measure on $ H / \Lambda $ clearly determines a finite
$J$ -invariant measure $\mu$ on the spaces $Y _z$ described above;
but these support no smooth invariant volume form. Indeed, suppose
$\nu$ were such a volume form, and let $ \frac {d \nu}{d \mu}$ be
the Radon-Nikodym derivative of (the measure defined by) $\nu$
with respect to $\mu$. Then this is a $J$-invariant function and
so, by ergodicity, constant almost everywhere; thus $\nu$ is a
constant multiple of $\mu$. But an easy calculation shows that $
\mu$ is (defined by) a smooth form on $Y _z$ which vanishes on the
submanifold $B$, since $\nu$ must be a constant multiple of this
form it is not a volume form.

However, imitating Katok-Lewis and Benveniste, we can modify the
action near $B$ to create a volume-preserving one. Let

$$\pi : J \times _{\Lambda_J} L \rightarrow J /\Lambda_J$$

\noindent be the natural projection, and $dH$ be the Haar volume
form on $H / \Lambda$. Then, considering $V - 0$ as a subset of
$L$, one finds easily that
$$
\phi ^* \mu = \pi ^* dH \wedge \Omega
$$
on $ W \cap H \times _{\Lambda} (V - 0)$ where $\Omega$ is a $J$ -
invariant volume form on the vector space $V$.  Now let $n$ be the
dimension of $V$. There is a $GL(V)$-equivariant (nonlinear) map
$k: L \rightarrow L ^{ \otimes n}$ defined by $v \rightarrow
\epsilon (v) ^ {(n+1)} v \otimes v \otimes ... \otimes v$, where
$\epsilon (v) = 1$ if $v$ has the same orientation as its image in
$S ^ +$, and $ = -1$ otherwise). Note that $L^{\otimes n}$ is a
line bundle and $k$ is a diffeomorphism on the complement of the
zero-section. Thus there is a map
\[
\displaystyle K : J \times _{\Lambda_J} L \rightarrow J \times
_{\Lambda_J} L ^ {\otimes n}
\]
which is a diffeomorphism on the complement of $J \times
_{\Lambda_J} S ^+$. Let $Y ^{\prime} _z$ be the manifold
\[
\displaystyle Y _z - B \sqcup W ^ {\prime \prime} / \equiv
\]
where $W ^ {\prime \prime}$ is a sufficiently small neigborhood of
the zero - section in $ J \times _{\Lambda_J} L ^ {\otimes n} $
and the relation $\equiv$ is generated by $w \equiv \phi \circ K ^
{-1} (w) $ for $w \in W ^ {\prime \prime} - J \times _{\Lambda_J}
S ^+$. Since the actions of $J$ correspond under this equivalence,
there is a smooth action $\rho ^{\prime} _z$ of $J$ on $Y
^{\prime} _z$.

Not that $Y'_z$ and $Y_z$ are easily seen to be diffeomorphic.  It
is also easy to see, using Ratner's theorem, that they are not $J$
equivariantly diffeomorphic, or even $F$ equivariantly
diffeomorphic when $F$ is any subgroup containing a unipotent
element. See \cite[Section 4]{FW} for more discussion of these
observations in the special case of the examples of Katok and
Lewis, the arguments there carry over more or less verbatim to
this case.

Now choose a positive definite inner product on $V$; this defines
trivializations of $L$ and all associated bundles, and a volume
$\Theta$ on $S ^+$; in terms of these trivializations, the map $k$
above can be written
$$
k: {\Ra} \times S ^+ \rightarrow {\Ra} \times S^+
$$
$$
(r, l) \rightarrow (sgn (r) ^{n+1} r ^n, l).
$$
Since $\Omega = c r ^{n-1} dr \wedge \Theta$, where $c$ is some
constant; we have $\Omega = k ^ * c ^{\prime} dr \wedge \Theta$ on
the complement of the zero section, and so $K ^ {-1 *} ( \pi ^* dH
\wedge \Omega)$ defines a form which extends to a nonzero volume
form on all of $W ^ {\prime \prime}$. Thus the form $\Sigma$ on $Y
^{\prime} _z$ defined by
\[
\displaystyle \Sigma \vert _{Y _z -B} = \mu
\]
\[
\displaystyle \Sigma \vert _{W ^ {\prime \prime}} = K ^ {-1 *} (
\pi ^* dH \wedge \Omega) \] defines a smooth volume on $Y
^{\prime} _z$ which is obviously $J$ invariant. Thus

\begin{proposition}
 The $J$ action $\rho'$ on $Y'_z$ preserves a smooth volume.
\end{proposition}

\section {Deformations of the exotic actions}
 \label{section:deformations}

We now wish to study how the actions constructed above depend on
the parameter $z$.  We fix a norm on $\fz$.  Let $z_0 \in Z$. For
$z \in Z$ close to $z_0$, $Y'_z$ is diffeomorphic to $Y'_{z_0}$,
which we abbreviate as $Y$. We can then consider that as $z$
varies we obtain varying actions

$$ \rho'_z : J \rightarrow \Diff(Y)$$

\noindent by choosing the diffeomorphisms of $Y'_z$ with $Y$
suitably, we can assume that $\rho'_z$ varies differentiably in
$z$. More explicitly, let $f$ be a function on $H/\Lambda$ that is
equal to one on a small neighborhood of $H[z_0]$ and 0 outside a
slightly larger neighborhood. Then for some $\varepsilon
> 0$, and any $v$ in $\fz$ with $\|v\|<\varepsilon$
the maps $\exp(fv)$ give a family of diffeomorphisms
 $$
\phi _v : H/ \Lambda - (J[z_0] \cup J[e]) \rightarrow H/ \Lambda -
(J[(\exp(fv)) z_0] \cup J[e]), \; \|v  \| < \epsilon
$$
These extend to diffeomorphisms
$$
\phi_v : Y'_{z_0} \to Y'_{(\exp(fv))z_0}
$$
such that if
$$
\alpha _z : J \to \Diff(Y'_z)
$$
are the $J$-actions defined above, we can define $\rho_v : J \to
Y'_{z_0}$ by
$$
\rho_v(j) = \phi_v ^{-1} \alpha_{z}(j) \phi_v
$$
\noindent where $z=\exp(v)z_0$.

The following is straightforward to check:

\begin{lemma}
The deformation $\rho_v$ is smooth; that is, the map
$$
B_{\fz}(0, \varepsilon) \times J \times Y \to Y
$$
given by $(v,h,y) \to (\rho_v (j) y)$ is smooth.  It is also
smooth when the actions are restricted to any closed subgroup of
$J$.
\end{lemma}

\noindent {\bf Remark:} We will also denote the action $\rho_v$ by
$\rho_z$ and $\phi_v$ by $\phi_z$ where $z=\exp(v)z_0$.

We now give an argument to prove Theorem \ref{theorem:simple} in
the case where $G$ admits a faithful, irreducible representation.
Note that this implies that the center of $G$ is cyclic, and so is
not true for all $G$.  The general case of Theorem
\ref{theorem:simple} will follow from the proof of Theorem
\ref{theorem:general} below. Let $G$ be our simple Lie group,
$\G<G$ any lattice and let $\pi:G{\rightarrow}GL(V)$ be a faithful
irreducible representation. Since $G$ is simple and connected,
$\pi$ takes values in $SL(V)$. We perform the construction above
with $m=\dim(V)$ so that $J=SL(V)=SL(m,\Ra)$ and $H=SL(n,\Ra)$
with $n>m+1$. We describe the centralizer of $G$ or $\G$ in $H$.

\begin{lemma}
\label{lemma:centralizer} The connected component of centralizer of
$G$ or $\G$ in $H$ is $DZ$ where $Z$ and $D$ are as above.
\end{lemma}

\begin{proof}
The Lie algebra $\fm$ of the centralizer of $G$ is just the
subalgebra of $\fh$ consisting of vectors which are invariant under
$\Ad_{H}|_{G}$ where $\Ad_H$ is the adjoint representation of $H$
into $GL(\fh)$.  By writing any element of $\fh$ as a block matrix
with diagonal blocks of size $n$ and $m$ respectively, it is to see
that $\Ad_H|_G$ splits as a direct sum of the trivial representation
on $\mathfrak{gl}(n-m)$, the restriction $\Ad_J|_G$ on
$\mathfrak{sl}(m)$, plus $n-m$ copies of $\pi$ plus $n-m$ copies of
the contragradient representation $\pi^*$.  Since $\pi$ is
irreducible, the set of invariant vectors is contained in
$\fz{\oplus}\fj {\oplus} \fd$.  All vectors in $\fz$ and $\fd$ are
invariant by definition, so it remains to see that there are no
invariant vectors in $\fj=\mathfrak{sl}(m)$.  The adjoint $\Ad_J$ is
just the restriction of the conjugation action of $J$ on
$\mathfrak{sl}(m)$ and $\mathfrak{sl}(m) \subset \mathfrak{gl}(m)$
which is isomorphic to $V{\otimes}V^*$ as a $J$ module.  Now
$V{\otimes}V^*$ is isomorphic to $\Hom(V,V)$ as a $J$ or $G$ module,
so by Schur's lemma and the fact that $\pi$ is irreducible, the set
of $G$ invariant vectors in $V{\otimes}V^*$ is exactly the scalars
and there are no $G$ invariant vectors in $\fj$ which is exactly the
complement of the scalars.

The centralizer of $\G$ is also an algebraic group whose Lie
algebra is the set of $\G$ invariant vectors in $\Ad_J|_{\G}$. By
the Borel density theorem, see e.g. \cite[Lemma $II.2.3$ and
Corollary $II.4.4$]{Ma}, this is the same as the set of $G$
invariant vectors, so we are done.
\end{proof}

We now state a more precise version of Theorem
\ref{theorem:simple}:

\begin{theorem}
\label{theorem:nontrivial} Given any $z_0$, the $G$ actions on
$Y_{z_0}$ and $Y_z$ are not conjugate as long as $zz_0{\inv}$ is
small enough. The same is true for the restriction of these actions
to $\G$.
\end{theorem}

\noindent {\bf Remark:} In order to show the existence of an
embedded disk in the deformation space, one needs to notice that the
bound required on the size of $zz_0{\inv}$ is continuous.  Therefore
by varying $z_0$ (as well as $z$) we can see that there is a
neighborhood of $\rho_{z_0}$ where no two actions are conjugate.

\begin{proof} Let $z \in Z$ be small enough so that $z,zz_0{\inv}{\notin}J\Lambda$, and
let $\psi: Y \to Y$ be a diffeomorphism conjugating the action
$\rho_0$ to the action $\rho_z$. By construction, the map $\phi_z$
conjugates the action $\rho_z$ to the action $\alpha_z$. Therefore
the map $\phi=\phi_z{\inv}{\circ}\psi$ conjugates the action
$\rho_0$ to the action $\alpha_z$.  After deleting the exceptional
divisor, the map $\phi$ is clearly an essentially injective,
essentially surjective map from $H/{\Lambda}$ to itself which
commutes with the action of $G$ or $\G$. By Corollary
\ref{corollary:rsw}, this implies that $\phi$ is covered by a map of
the form $A_{\phi}{\circ}L_{z_{\phi}}$ almost everywhere, where
$z_{\phi}$ is an element of $H$ and $A$ is an automorphism of $H$.
Since the only non-trivial automorphism of $H$, conjugate transpose,
does not commute with the action of $G$ or $\G$, we know that
$\phi=L_{z_{\phi}}$ almost everywhere where $z_{\phi}$ is an element
of $Z_H(G)=Z(H)DZ$. Since $\phi$ is smooth,
$\phi|_{Y\backslash{B}}=L_{z_{\phi}}$. We now show that $\phi$
cannot extend to a diffeomorphism of $Y$, obtaining a contradiction.

Since in our first copy of $Y$, we obtain $E$ by blowing up and
glued the orbits $J[e]$ and $J[z_0]$ and in our second copy of $Y$
we obtain $B$ by blowing up and gluing $J[e]$ and $J[z]$, for the
map $\phi$ to extend to $H/{\Lambda}$ would require that
$\phi(J[e])=J[e]$ and $\phi(J[z_0])=J[z]=zz_0{\inv}J[z_0]$ or
$\phi(J[e])=J[z]=zJ[e]$ and $\phi(J[z_0])=J[e]=z_0{\inv}J[z_0]$. In
either case, $z_{\phi}$ is in $Z(J)Z$ where $Z(J)$ is the center of
$J$. This is because $L_{z_{\phi}}$ maps $J[e]$ inside of $JZ[e]$
and therefore $z_{\phi}{\in}JZ$ by Corollary \ref{corollary:orbits}.
The intersection of $J$ with the centralizer of $G$ must lie in the
center of $J$ by Schur's lemma and therefore $z_{\phi}{\in}Z(J)Z$.
In the first case above, we have then have that
$L_{z_\phi}J[e]=J[z_{\phi}]=J[e]$ and
$L_{z_{\phi}}J[z_0]=J[z_{\phi}z_0]=J[z]$.  Since $JZ[e]$ is
equivariantly a product $J/{\Lambda_J} {\times} Z/{\Lambda_Z}$, the
first equation forces $z_{\phi}{\in}Z(J)\Lambda_Z$, while the second
forces $z_{\phi}=zz_0{\inv}$ modulo $Z(J)\Lambda_Z$ giving a
contradiction. The same analysis works in the second case.
\end{proof}

\noindent {\bf Remarks:} It is clear from the construction that we
have a space of deformations whose dimension is equal to the
dimension of $Z$.  Since $Z$ can be chosen to be arbitrarily
large, this provides examples with arbitrarily large deformation
spaces.  Similarly, by considering examples where we ``blow-up and
glue" along $k$ pairs of closed orbits, we can construct actions
with deformation spaces of dimension $k\dim(Z)$.

We now turn to the setting of Theorem \ref{theorem:general} where
$G$ is a general connected real algebraic group.  Recall from
subsection \ref{subsection:reductions} that $G=L{\ltimes}U$ where
$U$ is the unipotent radical and $L$ is reductive.  By
Propositions \ref{proposition:suffices} and
\ref{proposition:noncompact}, it suffices to consider the case
where $L$ is noncompact with compact center.  In this case, $L$ is
an almost direct product $Z(L)M$ where $Z(L)$ is compact and $M$
is semisimple with at least one noncompact factor.  We denote some
noncompact adjoint factor of $M$ by $G_0$ and note that there is a
homomorphism $\sigma:G{\rightarrow}G_0$.  We construct an action
of $G_0$ as above. First choose an irreducible representation of
$G_0$  which defines a map $\pi:G_0{\rightarrow}J=SL(m,\Ra)$ for
some $m$.  We then choose some $n>m+1$ and let $H_1=SL(n,\Ra)$. As
before the centralizer of $G_0$ in $H$ is an algebraic group with
connected component $Z=GL(n-m,\Ra)$ and we can construct a
cocompact lattice $\Lambda_1<H_1$ such that
$\Lambda_1{\cap}Z=\Lambda_1^Z$ is a lattice in $Z$,
$\Lambda_1{\cap}J=\Lambda_1^J$ is a lattice in $J$ and
$\Lambda_1{\cap}JZ=\Lambda_1^{JZ}$ is a lattice in $JZ$. The
construction above gives a family of actions $\rho_z^1$ of $G_0$
(and therefore of $G$) on a manifold $Y$ which is obtained from
$M_1=H_1/{\Lambda_1}$ by blowing up two closed $J$ orbits and
gluing.  Since $G$ is algebraic, we have an embedding of $G$ in to
$GL(N,\Ra)$ and therefore into $SL(N+1,\Ra)=H_2$.  We choose an
arbitrary cocompact lattice $\Lambda_2<H_2$.  Passing to a
subgroup of finite index in $\Lambda_2$ so that $Z(G)$ does not
intersect $\Lambda_2$, there is a faithful $G$ action $\rho^2$ on
$H_2/{\Lambda_2}$ by left translation.  Using $\rho_z^1$ and
$\rho^2$, we can define a diagonal action of $G$ on
$Y{\times}H_2/{\Lambda}$ which we will denote by $\rho_z$. We can
now state the following more precise version of Theorem
\ref{theorem:general}:

\begin{theorem}
\label{theorem:nontrivialgeneral} For a given $z_0$, the $G$
actions on $\rho_{z_0}$ and $\rho_z$ are not conjugate as long as
$zz_0{\inv}$ is small enough. The same is true for the restriction
of these actions to any lattice $\G<G$.
\end{theorem}

\begin{proof}
The proof is analogous to that of Theorem
\ref{theorem:nontrivial}.  The only difference is that applying
Corollary \ref{corollary:rsw}, we can only conclude that
$\phi=L_{z_{\phi}}$ where
$z_{\phi}{\in}Z_{H_1}(G_0){\times}Z_{H_2}(G)$. It is
straightforward to check that this does not effect our ability to
derive a contradiction.
\end{proof}

\section{More deformations of rank 1 lattice actions}
\label{section:extradeformations}

In this section, we describe some mixing actions of lattices in
$SO(1,n)$ which admit actions on homogeneous spaces for
$SO(1,n+1)$ for which there are non-trivial analytic deformations.
Note that since lattices in $SO(1,n)$ frequently admit
homomorphisms to $\Za$, it is trivial to construct non-faithful
actions which admit non-trivial perturbations following the proof
of Proposition \ref{proposition:suffices}. To proceed, we need to
recall a construction due to Johnson and Millson.  In \cite{JM},
they construct lattices $\G<SO(1,n)$ which admit continuous
families of non-trivial deformations when viewed as subgroups of
$SO(1,n+1)$. More precisely:

\begin{theorem}[Johnson-Millson]
\label{theorem:millsonjohnson} There exist both non-cocompact and
cocompact lattices $\G<SO(1,n)$ such that
\begin{enumerate}
\item $C=\G{\cap}SO(1,n-1)$ is a lattice in $SO(1,n-1)$ \item
there are subgroups $A,B<\G$ such that  $\G$ splits as a free
product with amalgamation $A*_{C}B$ \item let $Z$ be the
centralizer in $SO(1,n-1)$ in $SO(1,n+1)$, if $z_t{\in}Z$ is a
smooth path, then for $t$ small enough, the group $A*_C{B}^{z_t}$
is the image of a discrete faithful homomorphism
$\sigma_t:\G{\rightarrow}SO(1,n+1)$ whose image is not conjugate
to $\G$ in $SO(1,n+1)$.
\end{enumerate}
\end{theorem}

This construction is often referred to as bending the lattice. By
``bending" along more than one subgroup instead of just $C$,
Johnson and Millson construct lattices with large finite
dimensional deformation spaces.  We leave it to the reader to
verify that the construction below produces a deformation space
for the group action of the same dimension and only prove the
existence of non-trivial deformations.

If we choose a cocompact lattice $\Lambda<SO(1,n+1)$, the
homomorphisms $\sigma_t:\G{\rightarrow}SO(1,n+1)$ can be used to
define actions $\rho_t$ of $\G$ on $SO(1,n+1)/{\Lambda}$.  The
theorem of Johnson and Millson can now be interpreted as saying
that the actions $\rho_t$ are not conjugate in the group of affine
diffeomorphisms of $SO(1,n+1)/{\Lambda}$.  (It is easy to see that
the group of affine diffeomorphisms is just $SO(1,n+1)$.)  We now
show that the actions are actually not conjugate in the full
diffeomorphism group.  Examination of the first proof shows more:
the actions are not conjugate by any measurable isomorphism of
$SO(1,n+1)/{\Lambda}$.

\begin{theorem}
\label{theorem:bendingactions} The actions $\rho_{t_1}$ and
$\rho_{t_2}$ are conjugate if and only if $t_1=t_2$. Therefore,
for $t$ small enough, $\rho_t$ is a non-trivial deformation of
$\rho_0$.
\end{theorem}

{\noindent}{\bf Remarks:}
\begin{enumerate}
\item Both the theorem and it's proof also apply mutatis mutandi
to the other lattices for which Johnson and Millson produce
bending deformations.  These are lattices $\G$ where there is a
non-separating hypersurface rather than a separating hypersurface
and so $\G$ can be written as an $HNN$ extension, but not as a
free product with amalgamation.

\item I give two proofs below.  As mentioned in the
acknowledgements, the second proof given here was discovered
jointly with R.Spatzier.

\item A third proof of Theorem \ref{theorem:bendingactions}, by
computation of the entropy of the perturbed actions, was pointed
out to the author by Alex Furman. This proof, though it uses much
less technology, is somewhat more tedious, as it involves
computing entropy for elements not in either $A$ or $B$, and the
details are left to the interested reader.
\end{enumerate}

\begin{proof}[Proof $1$:] This proof applies only when $n>1$.
If $\rho_{t_1}$ and $\rho_{t_2}$ are conjugate, then there exists
$\phi$ in $\Diff(SO(1,n+1)/{\Lambda})$ such that
$\rho_{t_1}(\g){\circ}\phi=\phi{\circ}\rho_{t_2}$ for any
$\g{\in}\G$.  Therefore, for any $\g{\in}C$, we have that
$\rho_{t_1}(\g){\circ}\phi=\phi{\circ}\rho_{t_1}(\g)$ since
$\rho_{t_1}(\g)=\rho_{t_2}(\g)$ for $\g$ in $\G$. Then by Corollary
\ref{corollary:rsw} and the fact that $\phi$ is small and the
automorphism group of $SO(1,n+1)$ is discrete, we have that
$\phi=L_z$ where $z$ is in
$Z=Z_{SO(1,n+1)}(C)=Z_{SO(1,n+1)}(SO(1,n-1))$ since the centralizer
of $C$ or $SO(1,n-1)$ in the affine group is contained in the group
of translations. Then $\rho_2=T_{z{\inv}}{\circ}\rho_1{\circ}T_z$,
which implies that $\sigma_{t_1}=\sigma_{t_2}^{z}$.  Applying this
equation to any $\g{\in}A$, we see that ${\g}^z={\g}$, which implies
that $z$ is trivial since $A$ is Zariski dense by \cite[Lemma
5.9]{JM}. The same argument applied to $B$ implies that
$z_{t_2}z=z_{t_1}$, a contradiction unless $z_{t_1}=z_{t_2}$.
\end{proof}

\begin{proof}[Sketch of Proof $2$:]
This proof uses the existence of hyperbolic elements
$\g_1,{\ldots},\g_k$ in $A$ whose common centralizer in $SO(1,n)$
is trivial.  We leave it as an exercise for the reader to verify
this.  Each $\g_i$ has centralizer $\Ra_i{\times}SO(n-1)$ in
$SO(1,n)$ where $\Ra_i$ is the one parameter group containing
$\g_i$ and centralizer $\Ra_i{\times}{SO(n)}$ in $SO(1,n+1)$. The
intersection of these centralizers is just a copy of $SO(2)$. We
will use the theory of \cite{HPS} concerning normal hyperbolicity
and persistence of central foliations for partially hyperbolic
diffeomorphisms. Since the perturbed and unperturbed actions are
the same on $A$, we can denote the central foliation for the
(un)perturbed action of $\g_i$ by $W_i$.  The work of \cite{HPS}
now implies that any conjugacy $\phi$ between the perturbed and
unperturbed actions of $\G$ must map each leaf of each foliation
$W_i$ to itself, since the conjugacy is $C^0$ small and the $\g_i$
are normally hyperbolic with respect to these foliations for both
actions. This implies that $\phi$ must map each leaf of the
transverse intersection of these foliations to itself. This
transverse intersection is just the foliation $\ff$ defined by
$SO(2)$ orbits. Therefore $\phi$ can be written as $\phi(x)=z_xx$
where $z$ is in $SO(2)$. It is easy to check that the choice of
$z_x$ is invariant under $C$.  By ergodicity of the $C$ action, it
then follows that $z_x$ is constant, or that $\phi$ is translation
by $z$.  The proof finishes as before.
\end{proof}

\begin{corollary}
\label{corollary:so1n} There exist affine actions of $SO(1,n)$ on
homogeneous spaces of the form
$(SO(1,n){\times}SO(1,n+1))/({\Gamma{\times}\Lambda})$ with
nontrivial deformations.
\end{corollary}

\begin{proof} This follows from Theorem \ref{theorem:bendingactions} and Lemma \ref{lemma:induction}.  That
the induced action is an affine action on
$(SO(1,n){\times}SO(1,n+1))/({\Gamma{\times}\Lambda})$ can be
verified using the map
$$(SO(1,n){\times}(SO(1,n+1)/\Lambda)/\Gamma){\rightarrow}(SO(1,n){\times}SO(1,n+1))/({\Gamma{\times}\Lambda})$$
given by $[g,h]=[g,gh]$ which conjugates the induced action for
the unperturbed action to the diagonal action of $SO(1,n)$ on
$(SO(1,n){\times}SO(1,n+1))/({\Gamma{\times}\Lambda})$.
\end{proof}

\noindent{\bf Remark:} We note here that though the deformations
in Theorem \ref{theorem:bendingactions} are all affine, in
Corollary \ref{corollary:so1n}, only the unperturbed action is
affine.  More examples of bendings of actions of the groups
discussed in this section will be discussed in \cite{F}.

\noindent
Department of Mathematics\\
Indiana University\\
Rawles Hall\\
Bloomington, IN 47405\\


\begin{thebibliography}{99}

\bibitem{A}
Auslander, Louis On radicals of discrete subgroups of Lie groups.
Amer. J. Math. 85 1963 145--150.

\bibitem{B} E.J.Benveniste, Exotic Actions of Semisimple Groups and their
Deformations, unpublished manuscript and chapter of doctoral
dissertation, University of Chicago, 1996.

\bibitem{Be2}
E.J.Benveniste, Rigidity of isometric lattice actions on compact
Riemannian manifolds, {\it GAFA} 10 (2000) 516-542.

\bibitem{BF} E.J.Benveniste and D.Fisher, Nonexistence of
invariant rigid structures and invariant almost rigid structure,
{\it Comm. Anal. Geom.} 13 (2005) 89-111.

\bibitem{F} D.Fisher, Bending group actions and cohomology of
arithmetic groups, in preparation.

\bibitem{FM1}D.Fisher and G.A.Margulis, Local rigidity of affine actions of
higher rank groups and lattices, preprint.

\bibitem{FM2} D.Fisher and G.A.Margulis, Almost isometries, Property $(T)$ and
local rigidity, {\it Invent. Math.}, 162 (2005) 19-80. .

\bibitem{FM0}
D.Fisher and G.A.Margulis, Local rigidity for cocycles,in {\em
Surv.~Diff.~Geom.~Vol VIII}, refereed volume in honor of Calabi,
Lawson, Siu and Uhlenbeck , editor: S.T. Yau, 45 pages, 2003.

\bibitem{FW} Fisher, D., Whyte, K., Continuous quotients for lattice
  actions on compact spaces, {\it Geom. Dedicata}, 87 (2001), no. 1-3, 181--189.

\bibitem{HPS}
M.W. Hirsch, C.C. Pugh, and M. Shub, {\it Invariant Manifolds},
Lecture Notes in Mathematics 583 , Springer-Verlag, New York,
1977.

\bibitem{H} T.J.Hitchman, Rigidity theorems for large dynamical systems with
hyperbolic behavior. Ph.d. Thesis, University of Michigan, 2003.

\bibitem{JM} Johnson, Dennis; Millson, John J. Deformation spaces associated to
compact hyperbolic manifolds. Discrete groups in geometry and analysis
(New Haven, Conn., 1984), 48--106, Progr. Math., 67, Birkhäuser Boston, Boston, MA, 1987.

\bibitem{KL} Katok, A., Lewis, J.: Global rigidity
results for lattice actions on tori and new examples of
volume-preserving actions. Israel J. Math. 93 (1996), 253--280.

\bibitem{Ma} G.A. Margulis, {\it Discrete subgroups of semisimple Lie groups},
Springer-Verlag, New York, 1991.

\bibitem{MS} Myers, S. B.; Steenrod, N. E. The group of isometries of a
Riemannian manifold.{\it Ann. of Math.} 40 (1939), no. 2,
400--416.

\bibitem{Mo} Morris, Dave Witte, Introduction to Arithmetic Groups,
preprint, available at http://people.uleth.ca/~dave.morris/ or via
math arxiv.

%\bibitem{M}
%Mostow, G. D. Equivariant embeddings in Euclidean space. {\it Ann.
%of Math.} (2) 65 (1957), 432--446.

\bibitem{P} Palais, Richard S. Equivalence of nearby differentiable actions of a compact group.
{\it Bull. Amer. Math. Soc.}  67  1961 362--364.

%\bibitem{PS} Palais, Richard S.; Stewart, Thomas E. Deformations of compact
%differentiable transformation groups. {\it Amer. J. Math.} 82
%(1960) 935--937.


\bibitem{PR}
V. Platonov and A. Rapinchuk, {\it Algebraic Groups and Number
Theory}, Academic Press, New York, 1994.

\bibitem{R} M. Ratner, On Raghunathan's measure conjectures, {\it Ann. of
Math.} 134 no. 3 (1991) 545-607.

\bibitem{S} N. Shah, Invariant measures and orbit closures on homogeneous spaces
for actions of subgroups generated by unipotent elements, {\it
Proceedings of the International Colloquium on Lie Groups and
Ergodic Theory}
 Mumbai 1996, edited by S.G. Dani.


\bibitem{W} D. Witte, Measurable Quotients of Unipotent Translations on
Homogeneous Spaces, {\it  Trans. Amer. Math. Soc.} 354 no. 2
(1994) 577-594.

\bibitem{W2} D.Witte, Cocycle superrigidity for ergodic actions of non-semisimple
Lie groups. {\it Lie groups and ergodic theory (Mumbai, 1996)},
367--386, Tata Inst. Fund. Res. Stud. Math., 14, Tata Inst. Fund.
Res., Bombay, 1998.

\bibitem{ZB}
R.J. Zimmer,{\it Ergodic Theory and semisimple Groups},
Birkh\"auser, Boston, 1984.

\bibitem {Zim_ICM}  Zimmer, R.J.:Actions of semisimple groups and discrete
subgroups. Proc. Internat. Cong. Math. (Berkeley, 1986) 1247-1258 (1987)

\bibitem {Zim_Most}  Zimmer, R.J.:Lattices in semisimple groups and
invariant geometric structures on compact manifolds. In: Discrete groups in
geometry and analysis, R. Howe, editor. Prog. in Math {\bf 67},
152-210. Boston: Birkhauser 1987

\end{thebibliography}
\end{document}